\documentclass[10pt,twoside]{article}
\usepackage{amsfonts,epsfig}
\usepackage{amssymb}
\usepackage{times}

\newcommand{\cvd}{\hfill \mathrm{Q.E.D. }}
\newcommand{\qed}{\hfill Q.E.D. \\}

\newtheorem{theorem}{Theorem}[section]
\newtheorem{lemma}[theorem]{Lemma}

\newtheorem{proposition}[theorem]{Proposition}
\newtheorem{definition}[theorem]{Definition}

\newcommand{\reals}{\mathbb{R}}
\newcommand{\complex}{\mathbb{C}}

\newcommand{\module}[1]{\vert #1 \vert}
\newcommand{\Sdet}{\mathop{\mathrm{Sdet}}}
\newcommand{\pro}[2]{\langle #1,#2\rangle}
\newcommand{\diag}{\mathrm{diag}}

\begin{document}

\bibliographystyle{plain}
\title{
Left eigenvalues of $2\times 2$ symplectic matrices}

\author{
E.\ Mac\'{\i}as-Virg\'os
\and
M.J.\  Pereira-S\'aez}

\date{}

\pagestyle{myheadings}
\markboth{E.\ Mac\'{\i}as-Virg\'os and M.J.\  Pereira}{Left eigenvalues of $2\times 2$ symplectic matrices}
\maketitle

{\ }\thanks{Institute of Mathematics,
Department of Geometry and Topology, University of Santiago de Compostela, Spain
(quique.macias@usc.es). Partially supoorted by FEDER and   Research Project MTM2008-05861 MICINN Spain.}
\begin{abstract}
We obtain a complete characterization of the $2\times 2$ symplectic matrices having an infinite number of left eigenvalues.  Previously, we give a new proof of a result from Huang and So about the number of eigenvalues of a quaternionic matrix. This is achieved by applying an algorithm for the resolution of equations due to De Leo et al.
\end{abstract}

\\

\section{Introduction}
Left eigenvalues of $n\times n$ quaternionic matrices are still not well understood. For $n=2$, Huang and So gave in \cite{HuangSo-2001} a characterization of those matrices having an infinite number of left eigenvalues. Their result (see Theorem \ref{INFINITE} below) is based on previous explicit formulae by the same authors for solving some quadratic equations \cite{HuangSo-2002}. Later, De Leo et al. proposed in \cite{DeLeoDucatiLeonardi-2006} an alternative method of resolution which reduces the problem to finding the {\em right} eigenvalues of a matrix associated to the equation.

In this paper we firstly give a new proof of Huang-So's result,  based on the method by De Leo et al. Secondly, we completely characterize those symplectic matrices having an infinite number of left eigenvalues (see Theorem \ref{FINAL}). The application we have in mind is to compute in a simple way the so-called Lusternik-Schnirelmann category of the symplectic group $Sp(2)$. This will be done in a forthcoming paper \cite{GomezMaciasPereira-2008}.

\section{Left eigenvalues of quaternionic matrices}
Let 
$$A=\left[\matrix{a & b \cr c & d\cr}\right]$$ be a $2\times 2$ matrix with coefficients in the quaternion algebra $\mathbb{H}$. We shall always consider  $\mathbb{H}^2$ as a right vector space over $\mathbb{H}$.
\begin{definition}
{\rm A left eigenvalue of $A$ is any quaternion $q\in\mathbb{H}$ such that there exists some $u\in\mathbb{H}^2$, $u\neq 0$, with $Au=qu$.  }
\end{definition}

Clearly, if $bc=0$ then the eigenvalues are the diagonal entries. For a non-triangular matrix $A$, the following result  appears in \cite{HuangSo-2001} (with a different proof).

\begin{proposition} If $bc\neq 0$, the left eigenvalues of $A$ are given by $q=a+bp$, where $p$ is any solution of the unilateral quadratic equation \begin{equation}\label{EQUATION}
p^2+a_1p+a_0=0,
\end{equation}
with
$
a_1 = b^{-1}(a-d)$ and
$a_0 = -b^{-1}c$.
\end{proposition}

{\em Proof:\/}  Let us write $A=X+jY$, where $X,Y$ are complex $2\times 2$ matrices, and let 
$$\tilde{A}=\left[\matrix{ X & -\bar{Y}\cr Y & \bar{X}\cr}\right]$$ be the complexification of $A$. Let $\Sdet(A)=\module{\det(\tilde{A})}^{1/2}$ be the Study's determinant.
By using the axiomatic properties of $\Sdet$ (see \cite{Aslaksen-1996, CohenDeLeo-2000}) we can triangularize the matrix $A-qI$ obtaining that
$$\Sdet(A-qI) =\module{a-q}\cdot\module{
 (d-q)-c(a-q)^{-1}b}.
$$
It is easy to see that $q=a$ would imply $bc=0$. Hence $A-qI$ is not invertible iff
$(d-q)-c(a-q)^{-1}b=0$. By putting $p=b^{-1}(q-a)$ we obtain the desired equation.
 \qed 

In \cite{HuangSo-2001} Huang and So also proved the following theorem.
\begin{theorem} \label{INFINITE}
The matrix $A$ has either one, two or infinite left eigenvalues. The infinite case is equivalent to the conditions
$a_1,a_0\in\reals$, $a_0\neq 0$ and
$\Delta=a_1^2-4a_0<0$.
\end{theorem}

The original proof is based on a case by case study guided by the explicit formulae that the same authors obtained in \cite{HuangSo-2002} for solving an equation like (\ref{EQUATION}). In particular, they prove that in the infinite case the eigenvalues are given by te formula
$(a+d+b\xi)/2$, where $\xi$ runs over the quaternions $\xi\in\langle i,j,k \rangle$ with $\module{\xi}^2=\module{\Delta}$.

It is easy to see that the conditions in Theorem \ref{INFINITE} above are sufficient. In fact, 
if $a_0=s$, $a_1=t$, $s,t\in\mathbb{R}$,  then the definition of eigenvalue leads to  the equation
$q^2+tq+s=0$,
that, after the change $p=q+t/2$, gives  $p^2=\Delta/4<0$ which has infinite solutions $p=(\sqrt{-\Delta}/2)\omega$, $\omega\in\langle i,j,k\rangle$, $\module{\omega}=1$.

For the necessity of the conditions we shall give an alternative proof. It exploits an elegant method for the resolution of equations, proposed by De Leo et al. in \cite{DeLeoDucatiLeonardi-2006} as an improvement of a previous algorithm by Ser\^odio et al. \cite{SerodioPereiraVitoria-2001}

\section{The eigenvectors method}
In order to facilitate the understanding of this paper we explicitly discuss in this section the algorithm of De Leo et al. cited above.

Let $$M=\left[
\matrix{-a_1 & -a_0\cr 1 & 0 \cr}\right]$$
be the so-called {\em companion matrix} of equation (\ref{EQUATION}).
Then
$$M\left[\matrix{p\cr 1\cr}\right]=\left[\matrix{-a_1 p -a_0\cr p\cr}\right]=\left[\matrix{p^2\cr p\cr}\right]=\left[\matrix{p\cr 1\cr}\right]p,$$
which shows that in order to find the solutions we have to look for {\em right} eigenvalues $p$ of $M$ corresponding to eigenvectors of the precise form $(p,1)$. Accordingly to \cite{DeLeoDucatiLeonardi-2006} we shall call $p$ a {\em privileged} right eigenvalue.

\paragraph{Right eigenvalues}
The theory of right eigenvalues is well established \cite{Brenner-1951, Zhang-1997}. A crucial point is that the eigenvectors attached to a given right eigenvalue do not form a (right) $\mathbb{H}$-vector space.

\begin{proposition}\label{CONJUG}
  Let $\lambda$ be a right eigenvalue. Let $v$ be a $\lambda$-eigenvector. Then $vq^{-1}$ is a  $q\lambda q^{-1}$-eigenvector, for any  $q\in\mathbb{H}$, $q\neq 0$.
\end{proposition}

{\em Proof:\/}  Since $Mv=v\lambda$ we have $M(vq^{-1})=v\lambda q^{-1} =vq^{-1}(q\lambda q^{-1})$.
\qed 

As a consequence each eigenvector gives rise to a similarity class $[\lambda]=\{q\lambda q^{-1}\colon q\in\mathbb{H}, q\neq0\}$ of right eigenvalues. Recall that two quaternions $\lambda^\prime,\lambda$ are similar if and only if  they have the same norm, $\module{\lambda^\prime}=\module{\lambda}$, and the same real part, $\Re(\lambda^\prime)=\Re(\lambda)$. In particular, any quaternion $\lambda$ is similar to a complex number and to its conjugate $\bar \lambda$. 

\paragraph{Eigenvectors} So  in order to solve the equation (\ref{EQUATION}) we first need to find the {\em complex} right eigenvalues of the companion matrix $M$. These correspond to the eigenvalues of the complexified $4\times 4$ matrix $\tilde M$ and  can be computed by solving the characteristic equation $\det(\tilde M -\lambda I)=0$. Due to the structure of $\tilde M$  its eigenvalues appear in pairs $\lambda_1,\bar\lambda_1,\lambda_2,\bar\lambda_2$ \cite{Zhang-1997}. 

In order to compute the eigenvectors, let us consider the $\complex$-isomorphism 
\begin{equation}\label{TRASLATION}
(z^\prime,z)\in\complex^2 \mapsto z^\prime+jz\in\mathbb{H}.
\end{equation}

\begin{proposition} $(x^\prime,x,y^\prime,y)\in\complex^4$ is a $\lambda$-eigenvector of the complexified matrix $\tilde M$ if and only if $(x^\prime+jy^\prime ,x+jy)$ is a $\lambda$-eigenvector of $M$.
\end{proposition}

\paragraph{Equation solutions} Let $M$ be the companion matrix of equation (\ref{EQUATION}). Once we have found a complex right eigenvalue $\lambda$ of  $M$, and some $\lambda$-eigenvector $(q^\prime,q)$ we observe that $q^\prime=q\lambda$, due to the special form of $M$.  Hence by Proposition \ref{CONJUG} the vector
$$ \left[\matrix{q^\prime\cr q\cr}\right] q^{-1}= \left[\matrix{q \lambda q^{-1}\cr 1\cr}\right] $$
is a $q\lambda q ^{-1}$-eigenvector, that is $p=q \lambda  q^{-1}$  is a privileged eigenvalue in the similarity class $[\lambda]$ and hence it is the desired solution.

Notice that, by  Proposition \ref{CONJUG}, two $\mathbb{H}$-linearly dependent eigenvectors give rise to the same privileged eigenvalue.

\section{Number of solutions}
Now we are in a position to discuss the number of solutions of equation (\ref{EQUATION}). This will give  a new proof of Theorem \ref{INFINITE}.  

Let $V(\lambda)\subset \complex^4$ be the eigenspace associated to the eigenvalue $\lambda$ of the complexified matrix $\tilde M$. By examining the possible complex dimensions of the spaces $V(\lambda_k)$ and $V(\bar\lambda_k)$, $1\leq k \leq 2$, we see that:
\begin{enumerate}
\item
If the four eigenvalues $\lambda_1,\bar\lambda_1,\lambda_2,\bar\lambda_2$ are different, then each $V(\lambda_k)$ has dimension $1$ and gives  just one privileged eigenvalue $p_k$. Since $\bar\lambda_k$ gives $p_k$ too, it follows that there are exactly two solutions. 
\item
If some of the eigenvalues is real, say $\lambda_1\in\reals$, then all its  similar quaternions equal $p_1=\lambda_1$, independently of the dimension of $V(\lambda_1)$. So there are one or two solutions, depending on whether $\lambda_1=\lambda_2$ or not.
\item\label{TRES}
The only case where  infinite different privileged eigenvalues may  appear is when $\lambda_1=\lambda_2\not\in\reals$, which implies $\dim_\complex{V(\lambda_1)}=2=\dim_\complex{V(\bar\lambda_1)}$. 
\end{enumerate}

\paragraph{The infinite case}
So we focus on  case \ref{TRES}, when $\tilde M$ has exactly two different eigenvalues $\lambda_1,\bar\lambda_1$. 
The following Proposition proves that we actually have an infinite number of solutions (recall that $\lambda_1\notin\reals$).

\begin{proposition}\label{TOTAL}
 In case 3 all the quaternions similar to  $\lambda_1$ are privileged right eigenvalues of $M$.
\end{proposition}

{\em Proof:\/}  Take a $\complex$-basis $\tilde u,\tilde v$ of $V(\lambda_1)\subset\complex^4$ and the corresponding vectors $u,v$ in $\mathbb{H}^2$ by the isomorphism (\ref{TRASLATION}). Since the latter are of the form $(q\lambda,q)$, it follows that the second coordinates $u_2, v_2$ of $u$ and $v$ are $\complex$-independent in $\mathbb{H}$, hence a $\complex$-basis. That means that the privileged eigenvalues $p=q\lambda_1 q^{-1}$, where $q$ is a $\complex$-linear combination of $u_2$ and $v_2$, run over all possible quaternions similar to $\lambda_1$.
\qed

\paragraph{The Huang-So conditions} It remains to prove that in case 3 the conditions of Theorem \ref{INFINITE} are verified.

Let $\mathbb{H}_0\cong \mathbb{R}^3 $ be the real vector space of quaternions with null real part.  The scalar product is given by  $\langle q,q^\prime\rangle =-\Re(qq^\prime)$. An orthonormal basis is   $\langle i,j,k \rangle$. If $\xi\in\mathbb{H}_0$ we have  $\bar\xi=-\xi$ and $-\xi^2=\module{\xi}^2$.
Let $\Omega=S^3\cap \mathbb{H}_0$ be the set of vectors in $\mathbb{H}_0$ with norm $1$. It coincides with the quaternions similar to the imaginary unit $i$.

Let  $\lambda_1=x+iy$, $y\neq 0$, be one of the two complex eigenvalues of $\tilde M$ and let $p\in[\lambda_1]$ be any privileged eigenvalue of  $M$. Since $\Re(p)=\Re(\lambda_1)$ and $\module{p}=\module{\lambda_1}$, we can write $p=x+\module{y}\omega$, where $\omega$ is an arbitrary element of $\Omega$. 

Put  $a_1=t+\xi_1$, with $t\in\reals$ and $\xi_1\in\mathbb{H}_0$. From equation (\ref{EQUATION}) written in the form $a_0=-(p+a_1)p$  we deduce that
$$\Re(a_0)=xt+x^2-\module{y}^2+\module{y}\Re(\xi_1\omega).$$
Hence $\module{y}\pro{\xi_1}{\omega}$ does not depend on $\omega\in\Omega$. Since $y\neq 0$, the following Lemma \ref{RARO} ensures that  $\xi_1=0$, i.e. $a_1\in\reals$.

\begin{lemma}\label{RARO}
 Let $\xi\in \mathbb{H}_0$ verify $\pro{\xi}{\omega-\omega^\prime}=0$ for any pair $\omega,\omega^\prime$ of vectors in $\Omega$. Then $\xi=0$. 
\end{lemma}

{\em Proof:\/} 
Let $\xi=xi+yj+zk\neq 0$. Let us suppose for instance that $x\neq 0$ (the other cases are analogous). Take any $\omega\in \Omega$ orthogonal to $\xi$ and $\omega^\prime=i$. Then $\pro{\xi}{\omega-\omega^\prime}=x\neq 0$.
 \qed 

Now, a consequence of Proposition \ref{TOTAL} is that $p=\lambda_1$ is a privileged eigenvalue of $M$, so, by equation (\ref{EQUATION}), $a_0$ is a complex number. Since $\bar\lambda_1$ is a solution too, we deduce that $a_0=\bar a_0$ is a real number. Finally from $\lambda_1^2+a_1\lambda_1+a_0=0$ it follows that $a_1^2-4a_0<0$ because $\lambda_1\notin\reals$.

This ends the verification of the Huang-So conditions given in Theorem \ref{INFINITE}.

\section{Symplectic matrices}
Let us consider the $10$-dimensional Lie group $Sp(2)$ of $2\times 2$   symplectic matrices, that is quaternionic matrices
such that $A^*A=I$. Geometrically they correspond to the (right) $\mathbb{H}$-linear endomorphisms of $\mathbb{H}^2$ which preserve the hermitian product
\begin{equation}\label{HERMITIAN}
\langle u,v\rangle = u^*v =\bar u_1v_1+\bar u_2v_2.
\end{equation} 
Thus a matrix is symplectic iff their columns form an orthonormal basis for this hermitian product.

Let us find a general expression for any symplectic matrix.

\begin{proposition}\label{FORMA}
 A symplectic matrix $A\in Sp(2)$ either is diagonal or has the form
\begin{equation}
A=\left[
\matrix{\alpha &-\overline{\beta}\gamma\cr 
\beta &\beta\overline{\alpha}\overline{\beta}\gamma/ \vert \beta\vert^2}\right], \quad  \beta\neq 0, \vert \alpha\vert^2 + \vert \beta \vert^2=1, \vert \gamma \vert =1.
\end{equation}
\end{proposition} 

{\em Proof.}
From definition, the two columns $A_1, A_2$ of $A$ form an orthonormal basis of $\mathbb{H}^2$ for the hermitian product. Let the first one be
$$A_1=\left[ \matrix{\alpha\cr\beta\cr}\right], \quad \alpha,\beta \in \mathbb{H}, \quad \vert \alpha\vert^2 + \vert \beta \vert^2=1.$$
If
$\beta=0$, $A$ is a diagonal matrix 
$\diag(\alpha,\delta)$ with $\module{\alpha}=1=\module{\delta}$.

If $\beta \neq 0$, consider the (right) $\mathbb{H}$-linear map $\langle A_1,-\rangle \colon \mathbb{H}^2 \to \mathbb{H}$, which is onto because $\vert A_1\vert =1$. Then its kernel $K=(A_1)^\perp$ has dimension $\dim_\mathbb{H}K=1$.   Clearly, the vector
$$u=\left[\matrix{-\overline{\beta}\cr \beta\overline{\alpha}\overline{\beta}/ \vert \beta\vert^2}\right]\neq 0$$
where $\module{u}=1$, is orthogonal to $A_1$, so any other vector in   $K$ must be a quaternionic multiple of $u$. Since $A_2$ has norm $1$, we have
$$A_2 =u\gamma= \left[ \matrix{-\overline{\beta}\gamma\cr
\beta\overline{\alpha}\overline{\beta}\gamma/ \vert \beta\vert^2 } \right], \quad \gamma\in\mathbb{H}, \quad \vert \gamma \vert =1.~~~~~\cvd$$

\section{Left eigenvalues of a symplectic matrix} In this section we apply Theorem \ref{INFINITE} to the symplectic case. We begin with a
result that is also true for right eigenvalues.
\begin{proposition} The left eigenvalues of a symplectic matrix have norm $1$.
\end{proposition}

{\em Proof.} The usual euclidean norm in $\mathbb{H}^2\cong\reals^4$  is given by $\module{u}^2=u^*u=\langle u, u\rangle$, for the hermitian product $\langle,\rangle$ defined in (\ref{HERMITIAN}). Let $Au=qu$, $u\neq 0$. Then $$\module{u}^2=\langle u,u\rangle =  \langle Au, Au\rangle =\langle qu,qu\rangle =u^*\bar qqu =\module{q}^2\module{u}^2.~~~~~\cvd$$

Let $A\in Sp(2)$ be a symplectic $2\times2$ matrix. 
If $A=\diag(\alpha,\delta)$ is diagonal, then its left eigenvalues are $\alpha$ and $\delta$, as it follows from the equation
$$\Sdet(A-qI)=\vert \alpha -q\vert \cdot \vert \delta-q \vert =0.$$
Notice however that the right eigenvalues of $A$ are the similarity classes $[\alpha]$ and $[\delta]$,  which are infinite excepting when  $\alpha,\beta=\pm 1$.

\begin{theorem}\label{FINAL}
 The only symplectic matrices with an infinite number of left eigenvalues are those of the form
$$
\left[\matrix{ q\cos\theta & -q\sin\theta\cr
q\sin\theta &  q\cos\theta} \right],\quad   \quad\module{q}=1, \quad \sin\theta\neq 0.
$$
\end{theorem}

Such a matrix corresponds to the composition $L_q\circ R_\theta$ of a real rotation $R_\theta\neq\pm \mathrm{id}$ with  a left translation $L_q$, $\module{q}=1$.

{\em Proof:\/}  Clearly the matrix above verifies the Huang-So conditions of Theorem \ref{INFINITE}. Reciprocally, by taking into account Proposition \ref{FORMA} we must check those conditions for the values
$a=\alpha$, $b=-\bar\beta\gamma$, $c=\beta$, $d=\beta\bar\alpha\bar\beta\gamma/\module{\beta}^2$. Since 
$$a_0=-b^{-1}c=\bar\gamma\beta^2/\module{\beta}^2=s\in\reals$$ it follows that
$s=\module{\gamma}=1$ (notice that the condition $a_1^2-4a_0<0$ implies $a_0>0$). Then $$\gamma =\left(\beta / \module{\beta}\right)^2.$$
Substituting $\gamma$ we obtain $b=-\beta$ and $d=\beta\bar\alpha\beta/\module{\beta}^2$. Then
$a_0=-b^{-1}c =-\beta^{-1}\beta =1$.

We now compute
$$a_1=b^{-1}(a-d) = -\beta^{-1}(\alpha-\beta\bar\alpha\beta/\module{\beta}^2)={-1\over \module{\beta}^2}(\bar\beta\alpha-\bar\alpha\beta).$$
Hence $\Re(a_1)=0$, so the condition $a_1\in\mathbb{R}$ implies $a_1=0$. That means that $\bar\beta\alpha$ equals its conjugate $\bar\alpha\beta$, i.e. it is a real number. Call $r=\bar\beta\alpha\in\reals$. 

Since $\module{\alpha}^2+\module{\beta}^2=1$ and $\beta\neq 0$, it is $0<\module{\beta}\leq 1$. Take any angle $\theta$ such that $\module{\beta}= \sin\theta$, $\sin\theta\neq 0$. Define $q=\beta/\module{\beta}$, so we shall have $\beta=q\sin\theta$ with $\module{q}=1$. On the other hand,  the relationships $\module{\alpha}=\module{\cos\theta}$  and  $\module{r}=\module{\bar\beta}\module{\alpha}$ imply  that
$r=\pm \sin\theta\cos\theta$. By changing the angle  if necessary we can suppose that $r= \sin\theta\cos\theta$ without changing  $\sin\theta$. Then 
$$\alpha=r(\bar\beta)^{-1}=r\beta/\module{\beta}^2=q\cos\theta.$$

Finally $d= \beta\bar\alpha\beta/\module{\beta}^2=  q\cos\theta$ and the proof is done.
\qed 

Notice that the companion equation of a symplectic matrix is always $p^2+1=0$.








\end{document}